\documentclass{elsart}

\usepackage{amssymb,amsmath}
\usepackage{amsfonts}
\newtheorem{corollary}{Corollary}
\newtheorem{theorem}{Theorem}

\newtheorem{proposition}{Proposition}
\newtheorem{remark}{Remark}
\newtheorem{definition}{Definition}
\newcommand{\pp}{\noindent {\em Proof. }}
\newcommand{\bee}[1]{\begin{equation}\label{#1}}
\newcommand{\beq}[1]{\begin{eqnarray}\label{#1}}
\newcommand{\ene}{\end{equation}}
\newcommand{\eqe}{\end{eqnarray}}
\newcommand{\ld}{\ldots}

\newcommand{\cd}{\cdots}

\newcommand{\lb}[1]{\mathrm{Lab}(#1)}
\newcommand{\ra}{\rightarrow}

\newcommand{\wM}{\widetilde{M}}

\newcommand{\bx}{\hfill$\Box$}

\newcommand{\bM}{\overline{M}}

\newcommand{\far}[1]{\mathcal{A}\langle x_1,\ld,x_{#1}\rangle}
\newcommand{\rk}[1]{\mathrm{rank}\,#1}

\newcommand{\fd}{finite - dimensional}

\newcommand{\fa}{free associative algebra }

\newcommand{\fax}{\mathcal{A}\langle X\rangle }

\begin{document}
\begin{flushright}
\makebox[1.5in][r]{Dedicated to our teacher}\\
\makebox[2in][r]{\textbf{Alfred Lvovich Shmelkin}}\\
\makebox[1in][r]{\emph{on his 70th birthday}}
\end{flushright}
\begin{frontmatter}
\title{Schreier rewriting beyond the classical setting}

\author[yb,ybao,yb1]{Yuri Bahturin} and \author[ao,ybao,ao1]{Alexander Olshanskii}
\address[yb]{Department of Mathematics and Statistics, Memorial University of Newfoundland, St. John's, NL, A1C 5S7, \textsc{Canada}}
\address[ao]{Department of Mathematics,
1326 Stevenson Center,
Vanderbilt University, Nashville, TN 37240, \textsc{USA}}
\address[ybao]{Department of Algebra, Faculty of Mathematics and Mechanics, 119899 Moscow, \textsc{Russia}}
\thanks[yb1]{Partially supported by NSERC grant \# 227060-04 and URP grant, Memorial University of Newfoundland}
\thanks[ao1]{Partially supported by NSF grants DMS-0700811 and RFBR 08-01-00573}
\begin{abstract}
Using actions of the free monoids and free associative algebras, we establish some Schreier-type formulas involving the ranks of actions and the ranks of subactions in free actions or Grassmann-type relations for the ranks of intersections of subactions of free actions. The coset action of the free group is used to establish the generalization of the Schreier formula to the case of subgroups of infinite index. We also study and apply large modules over free associative algebras in the spirit of \cite{BO,OO}.
\end{abstract} 
\end{frontmatter}

\section{Introduction}

Let $\mathcal{A}_1,\ldots,\mathcal{A}_s$ be transformations of a set $M$. Then $M$ naturally acquires the structure of a $W$-set, where $W$ is a free monoid of rank $s$ or a free group of rank $s$, if these transformations are invertible or, finally,  a free associative algebra of rank $s$ (respectively, the group algebra of a free group of rank $s$) if $M$ is a vector space and all transformations are linear (respectively, linear and invertible).
 
Using actions of free monoids and free associative algebras, we establish some Schreier-type formulas involving the ranks of actions and the ranks of subactions in the free actions or Grassmann-type relations for the ranks of intersections of subactions of free actions. The coset action of the free group is used to establish a generalization of the Schreier formula to the case of subgroups of infinite index.

In the final portion of the paper we use Schreier techniques in the case of so called large modules over free associative algebras. The idea to use large objects, in particular groups, and Lie algebras has been successfully used in \cite{OO,BO} to produce examples for the problems of Burnside type. Here we use large modules to produce examples of finitely generated nil-modules with additional properties.

\section{Subacts of free acts over free monoids}\label{sSFA}

We start this section by recalling few known definitions and facts (see, e.g. \cite{MA}). By a \emph{monoid} we will understand a semigroup with identity element $1$.

\begin{definition} Given a monoid $S$, a nonempty set  $M$ is called a \emph{right act} over $S$ (or right $S$-act) if there is a map $\mu: M\times S\ra M$ satisfying the following conditions. If we write $\mu(m,s)=ms$ then 
\begin{enumerate}
	\item[\emph{(1)}] $m(st)=(ms)t$,
	\item[\emph{(2)}] $m=m1$.
\end{enumerate}
Here $m,s,t$ are arbitrary elements such that $m\in M$ and  $s,t\in S$.
\end{definition}

In this paper all $S$-acts  will be right. A subset $A$ is a \emph{generating set} of an $S$-act $M$ if $M=AS$. Any monoid $S$ is an act over itself if one chooses $\mu:S\times S\ra S$ to be the product in $S$. Given an $S$-act $F$ with a nonempty generating subset $A$ we say that $A$ is a \emph{basis} of $F$ if for any $a_1,a_2\in A$ and $s_1,s_2\in S$ it follows from $a_1s_1=a_2s_2$ that $a_1=a_2$ and $s_1=s_2$. An act $F$ possessing a (necessarily unique, up to permutation of elements!) basis $A$ is called \emph{free}. If $|A|=m$ then $F$ is \emph{free of rank $m$} and we write $\rk F=m$. By \cite[Theorem 5.13]{MA} any free act of rank $m$, $m$ finite or infinite, is isomorphic to the disjoint union of $m$ copies of free acts of rank 1, each isomorphic to the $S$-act $S$. 

The results we obtain in this section deal with free acts over free monoids. For any nonempty set (=alphabet) $X$, the set $W(X)$ of all words in $X$, including the empty word denoted by $1$, is a monoid under the juxtaposition (concatenation) $(u,v)\mapsto uv$. If $r$ is the cardinality of $X$ then one calls $W=W(X)$ the \emph{free monoid of rank $r$} and writes $\rk W=r$. One calls $X$ the \emph{basis} of $W$.

Given a subact $P$ of a free act $F$ over a free monoid $W$, there
is a simple Schreier-type formula (see Corollary \ref{cSPFAf}
below) for the rank of $P$ in terms of the ranks of $F$, $W$ and
the cardinality of $F\setminus P$, in the case where the values
$\rk F$, $\rk W$ and $|F\setminus P|$ are finite. But even in the
case where the above values are not necessarily finite we are able
to establish a similar relation (see Theorem \ref{tPSFA} below)
once we replace the numerical values by certain generating
functions called by some authors the \emph{Hilbert series}. This
approach can be traced back to \cite[\S 2.5]{PMC}, where P.M.Cohn
introduced the so called \emph{gosha} (after Golod - Shafarevich!)
of a right ideal of a ring with weak algorithm. Later this
approach was successfully used by V. Petrogradsky \cite{VP}, who
established Schreier-type relations in the case of subalgebras of
free Lie algebras and related algebras.

We start with a simple Nielsen-Schreier type result about the bases of subacts of free acts over free monoids. 

\begin{proposition}\label{pNSFA} Let $F$ be a free act with basis $A$ over a free monoid $W$ with a free generating set $X$ and a subact $P$. Then $P$ is free and a basis $B$ for $P$ consists of all elements $a\in A\cap P$ and $awx\in P$ (where $a \in A$, $w\in W$, $x\in X$) such that $aw\notin P$.
\end{proposition}

\pp Let $P'$ be a subact generated by $B$. If $P'\ne P$ then there is an element $au\in P$ with $u$ of minimal length, which is not in $P'$. If $u=1$ then $au\in A\cap P\subset B\subset P'$, a contradiction. Otherwise $u=u'x$, for some $x\in X$. If $au'\in P'$ then also $au=(au')x\in P'$, a contradiction. Otherwise, $au=(au')x\in B\subset P'$, a contradiction. So $B$ generates $P$. Now assume $aux, a'u'x'\in B$ and $auxv = a'u'x'v'$, for some $a,a'\in A$,$x,x'\in X$ and $u,v,u',v'\in W$. Then $a=a'$ and either $ux=u'x'$ or, say, $u'x'$ is a proper subword of $ux$ or $u'x'$ is a subword of $u$. But then $au\in P$, a contradiction. So $aux = a'u'x'$, hence $v=v'$, proving that $B$ is a basis of $P$.\bx

The construction of the basis $B$ of subact $P$ in the previous proposition is similar to the so called \emph{Schreier bases} in the case of subgroups of free groups or submodules of free modules over the free associative algebras \cite{JL}. This can be visualized if we view a free $W$-act $F$ with basis $A$ as a \emph{forest} (i.e. a disjoint union) of \emph{trees} $aW$, one for each $a\in A$. Each element $aw$ of $F$ is thus a vertex of $\mathcal{T}_a$ and, for each $x\in X$, there is a directed edge labeled by $x$ from $au$ to $aux$. We write $au\le aux$ and extend this relation by transitivity to the tree order on  $aW$ and a partial order on $F$. Using this identification, it becomes obvious that the basis of a subact $P$ of $F$ is just the set of all vertices of $P$ without predecessors in $P$.

Now let $M$ be a set endowed with a nonnegative integer valued \emph{degree function} such that for each $n=0,1,\ld$ the number $c_n(M)$ of elements $x\in M$ with $\deg x=n$ is finite. We define a formal power series (the \emph{Hilbert series}) $H(M,t)$ in one variable $t$ by setting
$$
H(M,t)=\sum_{n=0}^\infty c_n(M)t^n.
$$
The above series is finite if and only if $M$ is finite in which case $H(M,1)=|M|$, the number of elements in $M$. In particular, suppose there are degree functions as above on $W$ and $F$ such that $\deg (fv)=\deg(f)+\deg(v)$ and $\deg (uv)=\deg(u)+\deg(v)$ for any $f\in F$ and $u,v\in W$. Then by restriction they induce degree functions on $X$ and $A$. If $P$ is a subact in $F$, $B$ is the basis in $P$ then we also have  the degree functions on $B$ and $F\setminus P$. So we can speak about the well-defined Hilbert series $H(X,t)$, $H(A,t)$, $H(B,t)$, and $H(F\setminus P,t)$. 

\begin{theorem}\label{tPSFA} Let $P$ be a subact with basis $B$ of a free act $F$ with basis $A$ over a free monoid $W$ with basis $X$ such that their Hilbert series are well defined and $\deg x\ge 1$ for any $x\in X$. Then
\bee{ePSFA}
H(B,t)= H(A,t)+H(F\setminus P,t)(H(X,t)-1).
\end{equation}\bx
\end{theorem}

\pp We use the form of  the basis of a subact from Proposition \ref{pNSFA}. Then $c_0(F)=c_0(A)=c_0(A\cap P)+c_0(F\setminus P)$ and so $c_0(B) =c_0(A)-c_0(F\setminus P)$. If $n>0$ then the elements of $B$ of degree $n$, that is, the elements of $P$ of degree $n$  which have no predecessors, fall into two groups. Firstly, the elements of $A$ of degree $n$ which are not the elements of $F\setminus P$. Secondly, for each $i=0,1,\ld,n-1$, the elements of the sets $Y(n,i)$ consisting of the products $aux$ of degree $n$, such that $x\in X$, $\deg x=n-i$, and $aux\notin P$. Then we have the following formula.
\begin{eqnarray*}
c_n(B)&=& (c_n(A) - c_n(A\setminus(F\setminus P))) + \sum _{i=0}^{n-1} (c_i(F\setminus P)c_{n-i}(X) - c_n(Y(n,i)))\\ 
&=&c_n(A)+\sum_{i=0}^{n-1}c_i(F\setminus P)c_{n-i}(X)-c_n(F\setminus P). 
\end{eqnarray*}
If we rewrite this in terms of the Hilbert series, we obtain (\ref{ePSFA}).
\bx

The following is the precise analogue of the classical Schreier formula for subgroups of free groups, Lewin's formula for the submodules of free modules over free associative algebras \cite{JL}, Kukin's formula for restricted subalgebras of free restricted Lie algebras \cite{K}, and probably some other.

\begin{corollary}\label{cSPFAf} Let $P$ be a subact of a free act $F$ over a free monoid $W$ both having finite ranks, such that $F\setminus P$ is finite. Then also $P$ is of finite rank and
\begin{equation}\label{eSPFAf}
\rk P=\rk F+|F\setminus P|(\rk W-1).
\end{equation}
\end{corollary}
\pp We only need to notice that under the condition of the corollary, the Hilbert series $H(B,t)$ for the basis $B$ of $P$ is well defined and we have $H(B,1)=\rk P$.\bx

To conclude this section we would like to prove that the intersection of two finitely generated subacts of a free act is itself finitely generated. This result is similar to the famous Howson's Theorem in the theory of free groups. Another (nontrivial) case is free restricted Lie algebras \cite{K}. Actually, from our formula for the subacts of free acts over free monoids (\ref{ePSFA}) we easily derive the following formula, which is a full analogue of the Grassmann formula for the subspaces of vector spaces. In the following theorem, given a subact $P$ of a free act $F$ we write $B_P$ for the uniquely defined basis of $P$.

\begin{theorem}\label{tHTFA} Let $P$, $Q$ be subacts of a free act $F$ over a free monoid $W$ with basis $X$ such that the Hilbert series $H(P,t)$ and $H(Q,t)$ are well defined and $\deg x\ge 1$ for any $x\in X$. Then the Hilbert series of $B_P$, $B_Q$, $B_{P\cup Q}$, and $B_{P\cap Q}$ are well defined and 
\bee{e889}
H(B_P,t)+H(B_Q,t)=H(B_{P\cup Q},t)+H(B_{P\cap Q},t).
\ene
In particular, if $P$ and $Q$ are of finite ranks, then also $P\cap Q$ and $P\cup Q$ have finite ranks and
\bee{e890}
\rk P+\rk Q= \rk(P\cup Q)+\rk (P\cap Q).
\ene
\end{theorem}

\pp Let us use the set theoretic equation $A\setminus(A\cap B)=(A\cup B)\setminus B$. We can apply (\ref{ePSFA}) twice to obtain
$$
H(B_{P\cap Q},t) = H(P\setminus(P\cap Q),t)(H(X,t)-1) +H(B_P,t)
$$ 
and
$$
H(B_Q,t) = H((P\cup Q)\setminus Q,t)(H(X,t)-1) +H(B_{P\cup Q},t). 
$$
Subtracting one of these equations from the other and considering 
$$
H(P\setminus(P\cap Q),t)=H((P\cup Q)\setminus Q,t)
$$
we obtain the desired formula (\ref{e889}). If the ranks of $P$ and $Q$ are finite then both series of the left side of (\ref{e889}) are finite (positive!), hence both series on the right side are finite (positive!) and so both $\rk(P\cup Q)$ and $\rk(P\cap Q)$ are finite. Replacing $t$ by $1$ we obtain the desired formula (\ref{e890}). \bx

\section{Subgroups of infinite index in free groups}\label{sPSF}

A $G$-act $S$ over a group $G$ is customarily called a \emph{$G$-set} and a 1-generator subact $sG$ the orbit of $s\in S$. Given a subgroup $H$ of a group $G$, the set $G/H$ of right cosets of $H$ in $G$ is well-known to be a $G$-set via $(Hx)g=H(xg)$ for any $x,g\in G$. Any $G$-set is of course the disjoint union of orbits. If a subgroup $H$ is the stabilizer of $s\in S$ then the map $sg\rightarrow Hg$ is an isomorphism of $G$-sets $sG$ and $G/H$. Thus the free acts and their subacts that we studied in the previous sections are of little interest, since any free act is just the union of several copies of $G$ and any subact is the union of some of these copies. Of far greater interest are the $F$-sets $F/H$, where $H$ is a subgroup of the free group $F$. Studying them allows one to obtain important information about the subgroup $H$. One of the examples is the classical Schreier Formula for the rank of $H$ in terms of the rank of $F$ and the index $[F:H]$. If we assume $[F:H]$ infinite then, as in the previous section, we can form the Hilbert series $\mathcal{H}(F/H,t)$ and $\mathcal{H}(B,t)$, where $B$ is a free basis of $H$. Our formula (\ref{eGSFG}) below will link these two series together. Notice that the notion of degree will now be replaced by the length of reduced words in the free group. The length of a coset $Hg$ will be understood as the shortest length of reduced words representing the elements of $Hg$.  

In our study of $F/H$ in this section we will use the language of graphs.

We start with a free group  $F$ with a \emph{symmetric basis} $A$,
that is, a union $\{ a_1,\ld,a_r\}\cup\{ a_1^{-1},\ld,a_r^{-1}\}$,
where $\{ a_1,\ld,a_r\}$ is a free basis of $F$. Let $H$ be a
subgroup in $F$ of not necessarily finite index. Let us introduce
a directed coset graph with labelling
$\mathcal{G}=(V(\mathcal{G}),E(\mathcal{G}),\mathrm{Lab})$ as
follows. For vertices, we set $V=V(\mathcal{G})=F/H$. By the
previous paragraph, $V$ is the disjoint union of ``spheres'' $V_n$
defined as the set of cosets of length $n$. By definition, a  \emph{labelled edge} is a triple $e=(e_-,e_+,\lb{e})$, where the vertices $e_-$ and $e_+$ are called the \emph{source} and the \emph{target} of $e$, respectively, and $\lb{e}\in A$ the \emph{label} of $e$. In our case $E(\mathcal{G})$ consists of all triples $(Hg,Hga,a)$ where $Hg$ runs through $F/H$ and $a$ through $A$. So $e_{-}=f_{-}$ and $\lb{e}=\lb{f}$ always imply $e=f$. Two labelled edges $e$ and $f$ are called \emph{inverses}
of each other, $f=e^{-1}$, if $e_{+}=f_{-}$, $f_{+}=e_{-}$ and
$\lb{f}=\lb{e}^{-1}$.  If $e_{i}, e_{i+1}$ are two consecutive
edges in a path $p=e_1\cdots e_n$ then
$(e_{i+1})_{-}=(e_{i})_{+}=(e_{i}^{-1})_{-}$. If also
$\lb{e_{i+1}}=\lb{e_{i}}^{-1}=\lb{e_i^{-1}}$ then by definition
$e_{i+1}=e_{i}^{-1}$. As a result, if a path $p=e_1\cdots e_n$ is
reduced (that is, has no subpaths $ee^{-1}$) then its label
$\lb{p}=\lb{e_1}\cd\lb{e_n}$ is a reduced word in the free group
$F$.

Let us set $V(n)=\bigcup_{k=0}^nV_k$. We define $\mathcal{G}(n)$ as a complete subgraph of $\mathcal{G}$ with $V(\mathcal{G}(n))=V(n)$. Now let us construct a maximal subtree $\mathcal{T}$ as the union of $\mathcal{T}(n)$, where each $\mathcal{T}(n)$ is a maximal subtree in  $\mathcal{G}(n)$, $n=0,1,\ld$, and $\mathcal{T}(0)=\mathcal{G}(0)$. Once  $\mathcal{T}(n-1)$ constructed, we set $V(\mathcal{T}(n))=V(\mathcal{T}(n-1))\cup V_n$ and add to $E(\mathcal{T}(n-1))$ one (double) edge from $\mathcal{G}(n)$  connecting each vertex in $V_n$ with a vertex in $V_{n-1}$. Then $V(\mathcal{T})=V(\mathcal{G})$ and also the length of each coset $Hg$ coincides with the distance from $H$ to $Hg$, defined as the length of the unique reduced path $p$ from $H$ to $Hg$ in $\mathcal{T}$. We have that $\lb{p}\in Hg\in V_n$ is a reduced word of length $n$.  This allows us to define each``sphere'' $V_n$  as the set of vertices in the distance $n$ from $H$. Each set $V_n$ is finite and we can form the Hilbert series 
$$
\mathcal{H}(F/H,t)=\mathcal{H}(V,t)=\sum_{n=1}^\infty c_n(V)t^n\mbox{, where }c_n(V)=|V_n|.
$$

With $\mathcal{T}$ fixed, the set $E(\mathcal{G})$ splits into two subsets: the edges in $E(\mathcal{T})$, called \emph{tree} edges, and the edges in $E(\mathcal{G})\setminus E(\mathcal{T})$, called \emph{non-tree} edges. Suppose $e$ is a non-tree edge. Let $p$ be a path from $H$ to $e_{-}$ while $q$ a path from $H$ to $e_{+}$. Suppose $v=\lb{p}$, $a=\lb{e}$ and $w=\lb{q}$; then  $v a w^{-1}$ is called a \emph{Schreier generator} for $H$. Since $e^{-1}$ is also a non-tree edge, the Schreier generator defined by $e^{-1}$ is the inverse of the Schreier generator defined by $e$. The collection $B=B(\mathcal{T})$ of all Schreier generators is known to be a symmetric basis of $H$ \cite{LS}. As with the set of vertices, the symmetric basis $B$ can be written as the union of finite $B_n$, each consisting of elements of length $n$. Again we have
$$
\mathcal{H}(B,t)=\sum_{n=1}^\infty c_n(B)t^n\mbox{, where }c_n(B)=|B_n|.
$$

If $H$ is a subgroup of finite index $[F:H]$ in $F$ then the rank of $H$ can be expressed through $[F:H]$ and the rank of $F$ by means of the classical Schreier formula 
\bee{eCSF}
\rk H=(\rk F-1)[F:H]+1.
\ene 
In comparison with the cases of free acts and free modules over free ideal rings, the relation between the above two Hilbert series is more complex. An example we give at the end of this section even indicates that an explicit expression of the coefficients of the first series through the second one is not possible. Instead, one has to consider a formal series $\widetilde{\mathcal{H}}(B,t)=\sum_{n=1}^{\infty}a(n)t^n$, where, for each $n\ge 0$,  $$
a(n) = \frac{1}{4}\, c_{2n-2}(B)+\frac{1}{2}\, c_{2n-1}(B) + \frac{1}{4}\, c_{2n}(B)
$$
is the ``weighted'' sum of the values  representing the number of free generators (without ``doubling''!) of lengths $2n-2$, $2n-1$ and $2n$. Notice that in the case where the index $[F:H]$ is finite the sets $B$ and $F/H$ are finite, $\sum_{n=1}^\infty a(n)=\sum_{n=1}^\infty \frac{1}{2}\, c_n(B)$ and so replacing $t$ by $1$ in $\widetilde{\mathcal{H}}(B,t)$ and $\mathcal{H}(F/H,t)$ we obtain
\bee{values at 1}
\widetilde{\mathcal{H}}(B,1)=\frac{1}{2}\mathcal{H}(B,1)=\rk H\mbox{ and }\mathcal{H}(F/H,1)=[F:H].
\ene

With the above modifications in place, we can now prove the following generalization of the classical Schreier Formula.

\begin{theorem} Let $H$ be a subgroup of a free group $F$ of rank $r\ge 1$. Then there is a symmetric basis $B$ of $H$ such that 
\begin{equation}\label{eGSFG}
\widetilde{\mathcal{\mathcal{\mathcal{H}}}}(B,t)= \left(rt - \frac{t+1}{2}\right) \mathcal{\mathcal{H}}(F/H,t) +  \frac{t+1}{2}.
\end{equation}
\end{theorem}

\pp We set $b(n)=c_n(B)$,  $v(n)=c_n(V)$. Notice that if $n>0$ then from each vertex $Hg\in V_n$ there is exactly one edge of $\mathcal{T}$ going to $V_{n-1}$ (the inverse of the last edge in the path $p$ going from $H$ to $Hg$). Therefore, there are exactly $v(n)$ tree edges going from $V_{n-1}$ to $V_n$. Quite similarly, there are exactly $v(n+1)$ tree edges going from $V_n$ to $V_{n+1}$. 

With each non-tree edge of $\mathcal{G}$ with label $a\in A$ which goes from a vertex in $V_n$ to a vertex in $V_n$ one can associate a Schreier generator $u= va w^{-1}$  such that $|u|= 2n+1$, since the lengths of $v$ and $w$ equal $n$. Conversely, each Schreier generator $u$ of length $2n+1$ can be obtained in this manner. It follows that the number of such edges equals $b(2n+1)$.
 
Now with each non-tree edge $e$ from $V_n$ to $V_{n-1}$ one can associate a Schreier generator of length $2n$  (in this case $|v|=n$ and $|w|=n-1$). As a result, the number of these edges equals $b(2n)/2$  (an equal number of Schreier generators of length $2n$ is associated with the non-tree edges $e^{-1}$ going from a vertex in $V_{n-1}$ to a vertex in $V_n$).
 
Similarly, there are $b(2n+2)/2$ non-tree edges from $V_n$ to $V_{n+1}$.
 
The outcome of the previous argument is as follows. For each vertex $v\in V_n$ the number of edges $e$ for which $v=e_{-}$ is $2r$. Therefore, the total number of the edges $e$ such that $e_{-}=v$ for some $v\in V_n$ is $2rv(n)$. In this number, the contribution of tree edges is $v(n+1)+v(n)$, because $V_n$ is connected with $V_{n-1}$ and $V_{n+1}$ by $v(n)$ and $v(n+1)$ tree edges, respectively. The number of non-tree edges is $b(2n)/2 +b(2n+1) + b(2n+2)/2$. Finally, for any $n>0$, we have
$$
2rv(n) = v(n) +v(n+1)+ b(2n)/2 +b(2n+1) + b(2n+2)/2.
$$
If we recall the definition of $a(n+1)$, the previous equation takes the following form: 
\begin{equation}\label{(*)}
a(n+1) = ((2r-1)v(n) - v(n+1))/2.
\end{equation}

Now the following values are easy:
$$ 
v(-1)=0,\: v(0)=1,\: b(-2)=b(-1)=b(0)=0.
$$ 
It follows then that
\bee{(**)} 
a(0) = \frac{b(0)}{2} =0,\: a(1) = \frac{1}{2}\left(b(1) + \frac{b(2)}{2}\right).
\ene

Also notice that
\begin{equation}\label{(***)}
2r = b(1) + \frac{b(2)}{2} +v(1)=2a(1)+v(1).
\end{equation}
Indeed, the total number of edges from $H$ to $V_1$ is $2r$, the number of (non-tree) edges from $H$ to $H$ equals the number of Schreier generators $b(1)$, the number of tree edges from $H$ to $V_1$ equals the number $v(1)$ of length 1 elements of the Schreier transversal and the number of non-tree edges from $H$ to $V_1$ is the same as half the number of Schreier generators of length 2, that is, $b(2)/2$. 

To conclude the proof, we have to make sure the coefficients of each $t^n$ on both sides of (\ref{eGSFG}) are the same. For $t^n$ with $n>1$ this easily follows by (\ref{(*)}) while for $t^0$ and $t^1$ one has to use (\ref{(**)}) or (\ref{(***)}), respectively. \bx
 
In particular, if $[F:H]=j < \infty$ then considering (\ref{values at 1}), we obtain the classical Schreier formula (\ref{eCSF}).

An interesting particular case is that of so called ``even'' subgroups. We say that a subgroup $H$ of a free group $F$ is \emph{even} if $H$ is generated by elements of even length. Using our previous notation, this is equivalent to $b(2n+1) =0$, for all natural $n$. In this case $b(2n)+b(2n+2) =4a(n+1)$. We now set $\frac{1}{2}b(2n) = d(n)$ (this is the number of free generators of length $2n$) and form a series $\widehat{\mathcal{H}}(B,t) = \sum_{n=0}^{\infty} d(n) t^n$. We use $d(n)+d(n+1) = 2a(n+1)$ to obtain 
$$
(t+1)\widehat{\mathcal{H}}(B,t) =2\widetilde{\mathcal{H}}(B,t),\mbox{ hence }\widehat{\mathcal{H}}(B,t)  = \left(\frac{2rt}{t+1}-1\right) \mathcal{H}(F/H,t) +  1.
$$
(Consider that by (\ref{(**)}) we have $d(1) = \frac{1}{2}b(2) = 2r - v(1)$).
We could also write 
$$
\mathcal{H}(B,t)=2\left(\frac{2rt^2}{t^2+1}-1\right) \mathcal{H}(F/H,t^2) +  2.
$$

As a result, the Hilbert series for the set of Schreier generators in the case of even $H$ can completely be recovered if we know the Hilbert series for $F/H$.

Notice that in any case of $H$ - even or not - the recoverable values are  $$2a(n+1)  = b(2n)/2 +b(2n+1) + b(2n+2)/2,$$ each being the number of Schreier free generators in the ``double'' system of generators of the form $v a w^{-1}$, where the length of at least one of $v,w$ is $n$.

We conclude this section with the following.

\begin{remark}\label{r000} In the general case, if we know the grading by lengths on $F/H$, it is not possible to recover the grading by lengths on the set of Schreier generators for $H$.
\end{remark} 

To explain this, we consider a free group $F$ of rank $r$ with the symmetric set of generators $A$. Let $\mathcal{G}$ be a connected graph with $A$-labeling and a fixed vertex $v_0$ such that for each vertex $v$ and each label $a\in A$ there is one edge with $e_{-}=v$, $\lb{e}=a$ and one edge $f$ with  $f_{+}=v$, $\lb{f}=a$. Then $\mathcal{G}$ is the graph of cosets of a subgroup $H\subset F$ (the elements of $H$ can be read on the loops starting at $v_0$).

This allows us to give an  example, as follows. 

Let us assume that $H$ is a subgroup of the free group $F$ as above and such that in the corresponding graph $\mathcal{G}$ there is a non-tree edge with label $a$ from a vertex $o(1)\in V_{n-1}$ to a vertex $o(2)\in V_{n-1}$ and another non-tree edge, with the same label, from a vertex $o(3)\in V_n$ to a vertex $o(4)\in V_n$. Let us modify $\mathcal{G}$ by cutting out the above edges and pasting in two non-tree edges with the same label $a$, one from $o(1)$ to $o(4)$ and one more from $o(3)$ to $o(2)$. The resulting graph $\mathcal{G^{\prime}}$ satisfies the condition on the number of labeled edges incident to each vertex and if the original graph $\mathcal{G}$ was large enough, the new graph remains to be connected, as well. This gives rise to a subgroup $H^\prime$ with the same Schreier transversal as $H$ but with different sequence of numbers $b(2n-1)$, $b(2n)$ and $b(2n+1)$. Thus their Hilbert series of the systems of Schreier generators are different.

\section{Large modules over free associative and free group algebras}\label{sPSFM}

In this section we will be concerned with right $R$-modules, where $R$ is a free associative algebra or a free group algebra.  Recall that given a nonempty set $X$ and a field $\Phi$ one defines \emph{\fa} $R=\fax$ as the vector space over $\Phi$ whose basis is the free monoid $W=W(X)$. The product on $W(X)$ extends by distributivity to the whole of $R$. With free group $F(X)$ in place of $W(X)$ we obtain a \emph{free group algebra} $\mathcal{F}\langle X\rangle$. The cardinality of $X$ is called the rank of either of these algebras. Along with some free products of algebras, $\fax$ and $\mathcal{G}\langle X\rangle$ are examples of so called \emph{free ideal rings}, that is, rings $R$ in which every (right) ideal is a free $R$-module. A wealth of information about free ideal rings can be found in the monographs of P. M. Cohn \cite{PMC,PMCn} 

Speaking about the Schreier Formula, we have to note that in \cite{JL} Jacques Lewin had  established its exact analogue in the case of a free submodule $N$ of finite codimension in a finitely generated free module $M$ over a $\fax$ or  $\mathcal{F}\langle X\rangle$ with $|X|=r$:
\bee{eLSF}
\rk N=(r-1)\dim M/N+\rk M.
\ene

Note that as a corollary one can easily derive that a submodule of finite codimension in a finitely generated $K$-module is finitely generated, provided that $K$ is a finitely generated $\Phi$-algebra.

Later, in \cite[\S 2.5]{PMC} P. M. Cohn gave a formula for the generating function of the set of generators of an arbitrary right ideal of a ring with weak algorithm (these rings include free associative algebras). From this result, after some work, one can derive a Schreier Formula for the Hilbert series similar to (\ref{ePSFA}) in the case of submodules of free modules without the finiteness condition. In the statement of this result, whose proof we skip, the coefficient of $t^n$ in the Hilbert series of a subset $S$ of $F$ is defined as the number of elements of degree $n$ in $S$. 

\begin{theorem}\label{tPSFM} Let $R$ be a free associative algebra of rank $r$ over a field $\Phi$. Suppose $N$ is a submodule of a free (right) $R$-module $F$ of rank $s$ and $M=F/N$. Then for any free base $B$ of $N$  one has
$$
H(B,t) = H(M,t)(rt-1) +s.
$$\bx
\end{theorem}

With this abundance of information about submodules of free modules, we will restrict ourselves to the study of a wider class, called ``large modules''. In \cite{OO} and \cite{BO} the authors have considered large groups and large restricted Lie algebras. An amazing property of ``large'' objects is that they seem to be very close to the free ones and yet enable one to perform construction of infinite (-dimensional) objects with ``finiteness'' conditions, which makes them an appropriate tool for the solution of Burnside - type problems. In both these papers the Schreier-type arguments were an essential part. Following these  papers we give the following definition.

\begin{definition}\label{lm}
A right module $M$ over an algebra $R$ is called \emph{large} if $M$ contains a submodule $N$ of finite codimension, which admits a homomorphism onto $R$. 
\end{definition}

Now let $N^{\ast}=\mathrm{Hom}_R(N,R)$ be the (left) $R$-module of $R$-linear forms on $N$. In the case of free ideal rings, one can say that $M$ is large if and only if it has a submodule $N$ of finite codimension such that $N^{\ast}\ne\{ 0\}$. A module $N$ with $N^{\ast}=\{ 0\}$ is called \emph{bound} \cite{PMC}. In the case of a free ideal ring $R$ a module is bound if and only if it does not contain $R$ as a direct summand.

An example of a large yet bound module over a free associative algebra is as follows. Let $r,s$ be natural numbers, $r>1$, and suppose a module $M$ over a free associative algebra $R$ of rank $r$ is given by $M=\langle u_1,...,u_s\,|\,u_1x_1+\cdots+u_sx_s=0\rangle$. Then $M$ has a free submodule $N=\langle u_1x_1,...,u_1x_r, u_2,..,u_s\,|\,\rangle$ of codimension 1, hence large, but cannot be mapped onto $R$. Indeed, if $b_1,\ld, b_s$ are the images of $u_1,\ld, u_s$ in $R$ under such a homomorphism then $b_1x_1=-b_2x_2-\cdots -b_sx_s$ in $R$ so that $b_1=\ldots=b_s=0$. But then the homomorphism is not onto, a contradiction. 

\begin{remark} In the case of groups an example shows that a similarly defined large group does not necessarily admit a homomorphism onto a virtually free group, that is, a group that has a nontrivial free subgroup of finite index. We do not know if such example exists for restricted Lie algebras.
\end{remark}

\begin{definition}\label{dGDR}For any ring (algebra) $R$ we say that a right $R$-module $M$ is given by \emph{generators} $u_1$, $\ld$, $u_p$ and \emph{defining relations} $v_1=0,\ld,v_q=0$, $M=\langle u_1,\ld,u_p\,|\,v_1,\ld,v_q\rangle$ for short, if $M=F/N$, where $F$ is a free $R$-module with free basis $\{ u_1,\ld,u_p\}$ and $N$ a submodule generated by $v_1,\ld,v_q$.
\end{definition}

Actually, each of the sets: the set of generators and the set of relations can be infinite. If the set of generators is finite then $M$ is called \textit{finitely generated}, if the set of relations is finite then $M$ is called \textit{finitely related}. A finitely generated finitely related module is called \textit{finitely presented}.

Note that several results that follow (Proposition \ref{pHT}, Theorems \ref{tlm} and \ref{tLarge}) can be proven even in the case where the base field of coefficients $\Phi$ is not necessarily commutative (a skew-field). 

\begin{proposition}\label{pHT} Let $R$ be either a free associative algebra
or a free group algebra, both with free basis $\{ x_1,\ld,x_r\}$, $\Delta$ the augmentation ideal of $R$, $z_1,\ld,z_r$ the free generators of $\Delta$ as a free right $R$-modules, that is, $z_i=x_i$ in the first case and $z_i=x_i-1$ in the second, $i=1,\ld,r$. Then given a presentation of a module $M$ with $p$ generators and $q$ relations, there is an $s\ge 0$ and a presentation of $M$ of the type which we call \emph{affine}, with $p+s$ generators and $q+s$ relations, in which every \emph{relator} $v_i$ has the form 
$$
v_i=u_1f_{i1}+\ld+u_pf_{ip}
$$
where $f_{ij}$ is a polynomial in $z_1,\ld,z_r$ of degree at most $1$ for all $i=1,\ld,q$, $j=1,\ld,p$.
\end{proposition}

\pp In the case of free associative algebras we can apply a trick due to Higman which is described in \cite[5.8]{PMCn}. The same trick, when applied to a finitely presented module $M$ over a free group algebra $R=\Phi\left\langle x_1^{\pm 1},\ld,x_r^{\pm 1}\right\rangle$ produces a similar canonical form for the presentation of $M$ except that in the polynomials $f_{ij}$ we can have the entries of $x_1^{-1},\ld,x_r^{-1}$. These can be removed by introducing $pr$ new generators $u_{st}$ and new relations $u_s=u_{st}x_t$, for all $s=1,\ld,p$ and $t=1,\ld,r$. Obviously the new presentation will define the same module. Then we can replace in all ``old'' relations $v_1=0,
\ld,v_q=0$ each entry of $u_sx_t^{-1}$ by a new generator $u_{st}$. After we have done so, it is easy to rewrite all relations so that they become polynomials of degree at most $1$ in $z_1=x_1-1,\ld,z_r=x_r-1$. \bx
In the case of free ideal rings, the following result \cite{JL} holds.

\begin{theorem}\label{tfrm} If $R$ is a free ideal ring and $M$ a finitely related module then  $M$ has a free submodule of finite codimension. 
\end{theorem}

By this Theorem any finitely related module $M=F/N$ over a free associative or free group algebra $R$ of finite rank $r$ is either finite-dimensional or large. Moreover, if the number of relations is smaller than the number of generators then the module cannot be finite-dimensional. Indeed, if $\dim M=d<\infty$ then by Schreier - Lewin's formula (\ref{eLSF}), $\rk N=(r-1)d+\rk F\ge\rk F$. However, in applications it is important to have a more specific information about the factor-module of the large module $M$ by the submodule $Q$ which can be mapped onto a nonzero free module.
 
We start with the following technical result.

\begin{theorem}\label{tlm} Let $R$ be a free associative algebra or a free group algebra of rank $r>1$ and $\Delta$ the augmentation ideal of $R$. Suppose an $R$-module $M$ is given by an affine presentation with $p$ generators and $q$ relations such that $p-q>0$. Then M has a submodule $L$ which can be mapped onto $R$. One can choose $L$ of codimension at most $k$ with $M\Delta^k\subset L$, where $k=\left[\frac{q}{r}\right]+1$.
\end{theorem}

\pp Let $M$ be presented as in Definition \ref{dGDR}: $M=F/N$, where $F$ is free with free basis $\{ u_1,\ld,u_p\}$ and $N$ is free with free basis $\{ v_1,\ld,v_q\}$, $p-q>0$. Every relator has the form indicated in Proposition \ref{pHT}, where each $f_{ij}$ has degree at most 1. 
Let $T=[f_{ij}]$ be the \emph{relation matrix} of our presentation. Without loss of generality we may assume that the first columns of $T$ have only  constant entries. If these columns are linearly dependent, then by their elementary transformations we can make a column of zero. By switching to another basis of $F$ we then obtain a presentation for $M$, where one of the free generators does not enter defining relations. Obviously, then $M$ itself maps onto $R$. So we can assume the ``constant'' columns of $T$ linearly independent. We can use these and maybe some other columns to transform the columns of $T$ to produce a different column with constant term zero. This follows from $q<p$. We may assume this is the last column. This brings us to a new basis  $\{ u_1',\ld,u_p'\}$ of $F$ such that a certain number of the first generators enter the defining relations only with constant coefficients while $u_p'$ enters with coefficients having constant term zero. As a result,  all $v_i$ are in a submodule of $M$ generated by the elements
\begin{equation}\label{e001} 
u_1',\ld,u_{p-1}',u_p'z_1,\ld,u_p'z_r.
\end{equation}
The above elements, considered as the elements of $F$, form a free basis of a submodule $Q$. Now $M=F/N$ contains a submodule $M_1=Q/N$. We have $M\Delta\subset M_1$ and $\dim M/M_1\leq\dim F/Q=1$. If we could map $M_1$ onto a nonzero free $A$-module, our proof would be complete. Now $M_1$ is given by the generators (\ref{e001}) and the same relations $v_1=0,\ld,v_q=0$. The generators $u_p'z_1,\ld,u_p'z_r$ enter our relations with only constant coefficients, meaning that the number of constant columns in the matrix of relations of $M_1/N$ grew by $r$. If we repeat the same argument $k$ times so that $kr>q$, we will definitely have that $M_k/N$ maps onto a nonzero free module. The codimension of $M_k$ is at most $k=\left[\frac{q}{r}\right]+1$ and if we set $L=M_k$ we will obtain the desired submodule. \bx

Considering what was said in Proposition \ref{pHT}, we now easily obtain the following.
\begin{theorem}\label{tLarge}
Let $R$ be a free associative algebra or a free group algebra of rank $r>1$ and $\Delta$ the augmentation ideal of $R$. Suppose an $R$-module $M$ is given by a presentation with $p$ generators and $q$ relations such that $p-q>0$. Then $M$ is large. A submodule $L$ of finite codimension which can be mapped onto $R$ can be chosen in such a way that $M\Delta^k\subset L$ for almost all natural $k$.\bx
\end{theorem}

As an application of large modules, we will give a construction of infinite-dimensional finitely generated nil-modules over free associative algebras or free group algebras. We recall that a module $M$ over a ring $R$ with augmentation ideal $\Delta$ is called \emph{nil} if for any pair of elements $(x,u)\in M\times \Delta$ there is $n=n(x,u)$ such that $xu^n=0$. Such modules easily arise if we take an infinite-dimensional $r$-generator nil-algebra $A$ (see, e.g. \cite{G}) and view it as a module over $R=\far{r}$. However, one can show that the rate of growth of such modules drops significantly if compared to the growth of the free module (see our forthcoming paper \cite{BOMG}), whereas the growth of the module produced using our present method is very much the same as the growth of $R_R$.

 In the proof of our last theorem we use the following notation. Given two positive integers $k,m$ we introduce the homogeneous noncommutative polynomials of degree $m$ in $k$ variables $y_1,\ld,y_k$. Let $u_1,\ld,u_t$ be the set of commutative variables which also commute with $y_1,\ld,y_k$. For any sequence $(j_1,...,j_k)$ of nonnegative integers with $j_1+\cdots+j_k=m$, the polynomial $f_{j_1,\ld,j_k}(y_1,\ld,y_k)$ is defined as the coefficient of $u_1^{j_1}\cdots u_k^{j_k}$ in the expansion of $(u_1y_1+\cdots+u_ky_k)^m$. With $k$ fixed and $m$ growing, the growth of the number $d_k(m)$ of these polynomials is the same as the growth of the number of monomials of degree $m$ in $k$ commuting variables, that is, \emph{polynomial}. This should be compared with the growth of the number of monomials in $r$ variables of degree $m$, which is  $r^m$, hence \emph{exponential}, if $r>1$. Thus, for any $l\ge 0$, $k\ge 1$, and $r\ge 2$ there is $m$ such that $d_k(l+m)<r^m$.

Now we are ready to prove the following.

\begin{theorem}\label{twnm} Let $R$ be a free associative or a free group algebra of rank $r>1$, over a field $\Phi$. Then any $R$-module $M$ given by $p$ generators and $q$ defining relations with $p-q>0$ has a  homomorphic image $\bM$ which is an infinite-dimensional finitely generated nil-$R$-module. By construction, $\bM$ is the direct limit of large modules, with all structure homomorphisms surjective. Moreover, a factor-module $\wM$ of $\bM$ is still infinite-dimensional and nil but also residually finite-dimensional.
\end{theorem}

\pp We pick a (countable) linear bases $E=\{e_1, e_2,\ld\}$ in $M$ and   $B=\{b_1,b_2,\ld\}$ in the augmentation ideal $\Delta$. We then form the list $\mathcal{L}=\{\tau_1,\tau_2,\ld\}$ of all finite sequences $(e_i,b_{j_1},\ld,b_{j_k})$. We set $M_1=M$ and assume by induction on $t=1,2,\ld$, that, for $t\ge 1$, $M_t$ is a large $R$-module which is a factor-module of $M$ such that if $(e_i,b_{j_1},\ld,b_{j_k})$ is the $l^{\mathrm{th}}$ tuple then there is $n_l$ such that the image of $e_i(\lambda_1b_{j_1}+\cdots+\lambda_kb_{j_k})^{n_l}$ is zero in $M_t$ for all $l<t$, where $\lambda_1,\ld,\lambda_k$ are arbitrary elements of $\Phi$. Additionally, we assume that there is an integer $s(t)$ such that $M_t\Delta^{s(t)}$ admits a homomorphism onto $R$. If $t>1$ then we also have $s(t)>s(t-1)$. 

By Theorem \ref{tLarge} all these conditions are satisfied for $t=1$. To construct $M_{t+1}$ once $M_t$ has been defined, we consider the $t^{\mathrm{th}}$ tuple, which we denote by $\tau_t=(e,b_1,\ld,b_k)$.  We know that there is a surjective homomorphism $\alpha: M_t\Delta^{s(t)}\twoheadrightarrow R$. Then $\alpha(M_t\Delta^{s(t)+1})=\Delta$. Let us choose $m>1$ so that $d_k(s(t)+m+1)<r^m$. For any $q>0$, let $N_q$ stand for the submodule of $M_t$ generated by all 
$$
e f_{j_1,\ld,j_k}(b_1,\ld,b_k)\mbox{, where }j_1,\ld, j_k\in\{1,\ld,k\}\mbox{ and }j_1+\cdots+j_k=q.
$$ Then $\alpha(N_{s(t)+m+1})\subset\Delta^m$. By Theorem \ref{tLarge} we have that $P=\Delta^m/\alpha(N_{s(t)+m+1})$ has a submodule $Q$ of finite codimension which maps onto $R$ and such that for some $l>0$ we have $P\Delta^l\subset Q$. As a result we have that $P\Delta^l$ maps onto a submodule $T$ of finite codimension in $R$. Thanks to Schreier-Lewin formula (\ref{eLSF}) $T$ is necessarily a nonzero free $R$-module. Going back to $\Delta^m$, we find that $(\Delta^{l+m}+\alpha(N_{s(t)+m+1}))/\alpha(N_{s(t)+m+1})$ maps onto a nonzero free $R$-module. By the Isomorphism Theorem then $\Delta^{l+m}/\alpha(N_{s(t)+m+1})\cap\Delta^{l+m}$ maps onto free nonzero module.  Setting $n=s(t)+m+l$ and considering preimages, we determine that $M_t\Delta^{n}/N_{s(t)+m+1}\cap M_t\Delta^{n}$ maps onto $R$. Let us recall the $t^{\mathrm{th}}$ tuple $(e,b_1,\ld,b_k)$. Then for any $\lambda_1,\ld\lambda_k\in\Phi$ we have that $e(\lambda_1b_1+\ld+\lambda_kb_k)^n\in N_{s(t)+m+1}\cap M_t\Delta^{n}$. This allows one to set $M_{t+1}=M_{t}/N_{s(t)+m}\cap M_t\Delta^{n}$. In this module we have that $e(\lambda_1b_1+\ld+\lambda_kb_k)^n=0$ and also, if we set $s(t+1)=n>s(t)$ that a submodule $M_{t+1}\Delta^{s(t+1)}$ maps onto a free nonzero module.

To obtain a finitely generated infinite-dimensional nil-$R$-module which is the direct limit of large modules we define $\bM$ as the direct limit of 
$$
M_0\twoheadrightarrow M_1\twoheadrightarrow M_2\twoheadrightarrow\ld .
$$
Then, in this module, we have $e_i(\lambda_1b_{i_1}+\ld+\lambda_kb_{i_k})^n=0$, for any tuple, as described above, that is, $\bM$ is a nilmodule. On the other hand, let us show that $\bM$ is infinite-dimensional. Notice that all $M_t/M_t\Delta^{s(t)}$ are \fd\ and the chain $R\supset \Delta\supset\Delta^2\supset\ld$ is properly descending. Now, since $M_t\Delta^{s(t)}$ maps onto $R$ and $s(t+1) >s(t)$, it follows that $M_k\Delta^{s(t+1)}$ is a proper submodule in $M_t\Delta^{s_{t}}$. Therefore,
\begin{equation}\label{e002}
\dim M_{t+1}/M_{t+1}\Delta^{s(t+1)}=\dim M_{t}/M_{t}\Delta^{s(t+1)} > \dim M_{t}/M_{t}\Delta^{s_{t}},
\end{equation}
proving that, indeed, $\bM$ is  infinite-dimensional.

If we want to obtain a finitely generated infinite residually finite nil-module, we have to consider 
$$
\widetilde{M}=\bM\left/\bigcap_{k=0}^\infty \bM\Delta^{s(t)}\right..
$$
It is obvious that $\wM$ is nil and residually finite-dimensional. Also notice that the kernel of the natural map of $\wM$ onto $M_t$ is contained in $\wM\Delta^{s(t)}$, for every $k$. Hence $\wM/\wM\Delta^{s(t)}\cong M_t/M_t\Delta^{s(t)}$. Now, by (\ref{e002}), $\dim M_t/M_t\Delta^{s(t)}\rightarrow\infty$. Since $\wM$ maps homomorphically onto every $M_t/M_t\Delta^{s(t)}$, it follows that $\wM$ is infinite-dimensional.
\bx

\end{document}